\documentclass[reqno,11pt]{article}
\usepackage{amsmath,amssymb}
\usepackage{amsthm,amscd,amsfonts,}
\usepackage{amsmath}
\usepackage{amsthm}
\usepackage{graphicx}
\usepackage{amsmath,amsthm}
\usepackage{amsmath,amstext,amsfonts}
\usepackage{subfigure}
\usepackage[T1]{fontenc}
\usepackage[utf8]{inputenc}
\usepackage{authblk}

\newtheorem{theorem}{Theorem}[section]
\newtheorem{lemma}{Lemma}[section]
\newtheorem{remark}{Remark}[section]
\newtheorem{corollary}{Corollary}[section]
\newtheorem{example}{Example}[section]
\newtheorem{definition}{Definition}[section]

\newcommand{\R}{\mathbb{R}}
\numberwithin{equation}{section}

\begin{document}

\title{\textbf{ Exactness of Second Order Ordinary Differential Equations and Integrating Factors}}
\author{R. AlAhmad\thanks{rami\_thenat@yu.edu.jo}}
\author{M. Al-Jararha\thanks{mohammad.ja@yu.edu.jo}}
\author{H. Almefleh\thanks{almefleh@yu.edu.jo }}
\affil{Department of Mathematics, Yarmouk University, Irbid, Jordan, 21163.}

\date{}
\maketitle

\begin{abstract}
A new problem is studied, the concept of exactness of a second order nonlinear  ordinary differential equations is established. A method is constructed to reduce this class into a first order equations. If the second order equation is not exact we introduce, under certain conditions, an integrating factor that transform it to an exact one.
\end{abstract}
 

\vspace{0.2in} \noindent \textbf{Keywords and Phrases}: Second order differential equation,
Exact equations, Non-exact equations,  Integrating factor.

\noindent \textbf{AMS (2000) Subject Classification}: 65L05, 49K15, 37C10


\section{Introduction}
The concept of exactness for a class of a first order nonlinear differential equations was presented \cite{Jordan} with a well-defined method of solution. The notion of integrating factor were introduced to convert differential equation that is not exact into an exact one. 

Second order nonlinear differential equations play an important role in Applied Mathematics, Physics, and Engineering \cite{Ames,Davis, Ghandraskhar,Hale,Jordan, Laselle,Lefschetz,Schlichting,Struble,von}. To find the general solution of a nonlinear  second order differential equation is not an easy problem in the general case. In fact, a very specific class of nonlinear second order differential equations can be solved by using special transformations. Another approach to study the solution of nonlinear second order differential equations is the dynamical systems approach. Using this approach a qualitative solution is given  instead of the particular solution of the equation.  A class of these equations will be solved in this paper. 

The outline of the paper: we give mathematical formulation for the exactness of a class of second order nonlinear equations based on transforming them into a first order equations. Also, we will introduce the idea of integrating factor to convert some differential equations into exact equations, and we will prove some related results.


\section{Exact Second Order Differential Equations}
Consider the following nonlinear second order differential equation
\begin{equation}\label{nonlinearsecondorder}
a_2(x,y,y^\prime)y^{\prime\prime}+a_1(x,y,y^\prime)y^\prime+a_0(x,y,y^\prime)=0.
\end{equation}
If a function $\Psi(x,y,y^\prime)$ exists with the properties that
\begin{equation}\label{exactcondition}
\frac{\partial \Psi(x,y,y^\prime)}{\partial x}=a_0(x,y,y^\prime), \; \frac{\partial \Psi(x,y,y^\prime)}{\partial y}=a_1(x,y,y^\prime), \;
\text{and}\;\; \frac{\partial \Psi(x,y,y^\prime)}{\partial y^\prime}=a_2(x,y,y^\prime),
\end{equation}
then we have 
\begin{equation*}
\frac{\partial \Psi(x,y,y^\prime)}{\partial y^\prime}y^{\prime\prime}+\frac{\partial \Psi(x,y,y^\prime)}{\partial y}y^\prime+\frac{\partial \Psi(x,y,y^\prime)}{\partial x}=0.
\end{equation*}
So, by the chain rule, we get
 $$\frac{d \Psi(x,y,y^\prime)}{d x}=0.$$
Hence,
\[
 \Psi(x,y,y^\prime)=c
\]
reduces Eq. \eqref{nonlinearsecondorder} into a first order differential equation.
 
\begin{definition}
The nonlinear second order differential equation \eqref{nonlinearsecondorder} is called exact equation  if there exists a function $\Psi(x,y,y^\prime)$ such that \eqref{exactcondition} holds.
\end{definition}

\begin{theorem}
Let the functions  $a_2(x,y,y^\prime)$, $a_1(x,y,y^\prime)$, $a_0(x,y,y^\prime)$,
$\displaystyle\frac{\partial a_2}{\partial y}$,  $\displaystyle\frac{\partial a_2}{\partial x}$, $\displaystyle\frac{\partial a_1}{\partial x}$, $\displaystyle\frac{\partial a_1}{\partial y^\prime}$, $\displaystyle\frac{\partial a_0}{\partial y^\prime},$ and $\displaystyle\frac{\partial a_0}{\partial y}$ be continuous functions in a simply connected region $R\subseteq \R^3$. Then Eq. \eqref{nonlinearsecondorder} is exact if and only if
\begin{equation}\label{exactcond2}
\frac{\partial a_2}{\partial y}=\frac{\partial a_1}{\partial y^\prime}, \;\; \frac{\partial a_2}{\partial x}=\frac{\partial a_0}{\partial y^\prime},\; \text{and}\; \;\frac{\partial a_1}{\partial x}=\frac{\partial a_0}{\partial y}.
\end{equation}

\end{theorem}
\noindent \textit{Proof.}  Assume that \eqref{exactcond2} hold. We are going to construct a function $\Psi(x,y,y^\prime)$ such that  
\[
\frac{\partial \Psi(x,y,y^\prime)}{\partial x}=a_0(x,y,y^\prime).
\]
Then by integrating this equation with respect to $x$, we get
\begin{equation}\label{proofequ}
\Psi(x,y,y^\prime)=\int_{x_0}^{x}a_0(\alpha,y,y^\prime)d\alpha+\Phi(y,y^\prime).
\end{equation}
Therefore, by differentiating the above equation with respect to  $y$ and using the assumption, we get
\[
\frac{\partial \Phi(y,y^\prime)}{\partial y}=a_1(x_0,y,y^\prime).
\]
Hence, 
\[
\Phi(y,y^\prime)=\int_{y_0}^{y}a_1(x_0,\beta,y^\prime)d\beta+\xi(y^\prime).
\]
To find $\xi(y^\prime)$, we substitute $\Phi(y,y^\prime)$ in Eq. \eqref{proofequ} to get

\[
\Psi(x,y,y^\prime)=\int_{x_0}^{x}a_0(\alpha,y,y^\prime)d\alpha+\int_{y_0}^{y}a_1(x_0,\beta,y^\prime)d\beta+\xi(y^\prime).
\]
Differentiate this equation with respect to $y^\prime$ and again use the assumptions, to get
\[
\xi^\prime(y^\prime)=a_2(x_0,y_0,y^\prime).
\]
Therefore, 
\[
\xi(y^\prime)=\int_{y^\prime_0}^{y^\prime}a_2(x_0,y_0,\gamma)d\gamma.
\]
Hence,
\[
\Psi(x,y,y^\prime)=\int_{x_0}^{x}a_0(\alpha,y,y^\prime)d\alpha+\int_{y_0}^{y}a_1(x_0,\beta,y^\prime)d\beta+\int_{y^\prime_0}^{y^\prime}a_2(x_0,y_0,\gamma)d\gamma.
\]
 The proof of the other direction is obvious. In fact, it comes from the assumption that 
$a_2(x,y,y^\prime)$, $a_1(x,y,y^\prime)$, and $a_0(x,y,y^\prime)$ are continuous with their first partial derivatives.
\begin{remark}
From the above theorem, we conclude that the nonlinear second order differential equation \eqref{nonlinearsecondorder} is  exact equation  if the conditions
\begin{equation}\label{exactconditions}
\frac{\partial a_2}{\partial y}=\frac{\partial a_1}{\partial y^\prime}, \;\; \frac{\partial a_2}{\partial x}=\frac{\partial a_0}{\partial y^\prime},\; \text{and}\; \;\frac{\partial a_1}{\partial x}=\frac{\partial a_0}{\partial y}
\end{equation}
hold.
\end{remark}
\begin{example}(The Plane Hydrodynamic Jet)
Consider the  second order nonlinear differential equation
$$ 3\epsilon y^{\prime\prime}+yy^\prime  =0.$$ This is exact. By using the result in the above theorem, we have
\begin{eqnarray}
\Psi(x,y,y^\prime)&=&\int_{0}^{y}\beta d\beta+3\epsilon\int_{_0}^{y^\prime}d\gamma\\
&=&\frac{y^2}{2}+3\epsilon y^\prime.
\end{eqnarray}
Hence, the equation is reduced to $\Psi(x,y,y^\prime)=c^2$, which is equivalent to 
\[
3\epsilon y^\prime+\frac{y^2}{2}=c^2.
\]
\end{example}
\begin{remark}
Consider the following nonlinear second order differential equation 
\begin{equation}\label{example1equation}
y^{\prime\prime}+a_1(x,y)y^\prime+a_0(x,y)=0,
\end{equation}
where $a_1(x,y)$ and $a_0(x,y)$ satisfy the condition \eqref{exactcondition}. Note that $a_2(x,y,y^\prime)=1$, $a_1(x,y,y^\prime)=a_1(x,y)$, and $a_0(x,y,y^\prime)=a_0(x,y)$, and so, it is obvious to see that the conditions \eqref{exactconditions} hold. Therefore, \eqref{example1equation} is exact.
\end{remark}

\begin{example} The  second order nonlinear initial value problem 
\begin{equation}\label{exampleeq}
\left\{
\begin{array}{l}
y^{\prime\prime}+12xy^3y^\prime+\left(3y^4-1\right)=0\\
y(0)=2, \;\;y^\prime(0)=0,
\end{array}
\right.
\end{equation}
is exact. Therefore, there exists a  function $\Psi(x,y,y^\prime)$ which reduces the above equation into a first order differential equation. By applying the above theorem, we have 
\[
\Psi(x,y,y^\prime)=\int_{x_0}^{x}a_0(\alpha,y,y^\prime)d\alpha+\int_{y_0}^{y}a_1(x_0,\beta,y^\prime)d\beta+\int_{y^\prime_0}^{y^\prime}a_2(x_0,y_0,\gamma)d\gamma,
\]
and since $x_0=0$, $y_0=2$, and $y_0^\prime=0$, we have
\begin{eqnarray*}
\Psi(x,y,y^\prime)&=&\int_{0}^{x}a_0(\alpha,y,y^\prime)d\alpha+\int_{2}^{y}a_1(0,\beta,y^\prime)d\beta+\int_{0}^{y^\prime}a_2(0,2,\gamma)d\gamma,\\
&=&\int_{0}^{x}(3y^4-1)d\alpha+\int_{0}^{y^\prime}d\gamma,\\
&=& y^\prime+(3y^4-1)x.\\
\end{eqnarray*}
Hence,  $\Psi(x,y,y^\prime)=c$ reduces  Eq. \eqref{exampleeq} to 
\[
y^\prime+(3y^4-1)x=c.
\]
By applying the initial data, we get $c=0$. Hence, Eq. \eqref{exampleeq} is reduced to the following first order differential equation
\[
y^\prime+3xy^4-x=0.
\]
For which an implicit solution can be obtained by separating the variable. 
\end{example}

%
\section{Non-exact Second Order Differential Equations and  Integrating Factors}
In this section, we introduce the notion of the integrating factor for the second order differential equation \eqref{nonlinearsecondorder}. Also, we deduce some conditions for the existence of such integrating factor. First, we start by the following definition for the integrating factor:
\begin{definition}
An integrating factor of Eq. \eqref{nonlinearsecondorder} is a non zero function $\mu(x,y,y^\prime)$, such that the equation
\begin{equation}\label{munonlinearsecondorder}
\mu(x,y,y^\prime)a_2(x,y,y^\prime)y^{\prime\prime}+\mu(x,y,y^\prime)a_1(x,y,y^\prime)y^\prime+\mu(x,y,y^\prime)a_0(x,y,y^\prime)=0
\end{equation}
is exact. i.e., 
\begin{equation}\label{muexactconditions}
\frac{\partial A_2}{\partial y}=\frac{\partial A_1}{\partial y^\prime}, \;\; \frac{\partial A_2}{\partial x}=\frac{\partial A_0}{\partial y^\prime},\; \text{and}\; \;\frac{\partial A_1}{\partial x}=\frac{\partial A_0}{\partial y},
\end{equation}
where
\[
A_2(x,y,y^\prime)=\mu(x,y,y^\prime)a_2(x,y,y^\prime),
\]
\[
A_1(x,y,y^\prime)=\mu(x,y,y^\prime)a_1(x,y,y^\prime),
\]
and
\[
 A_0(x,y,y^\prime)=\mu(x,y,y^\prime)a_0(x,y,y^\prime).
\]
\end{definition}
\begin{theorem}\label{thm31}
Assume that Eq. \eqref{nonlinearsecondorder} is not an exact equation. Then, it 
has no integrating factor of one of the forms $\mu(x,y,y^\prime)$, $\mu(x,y)$, $\mu(x,y^\prime)$,  or $\mu(y,y^\prime)$ if and only if
\begin{equation}\label{condformu}
\left(\frac{\partial a_0}{\partial y}- \frac{\partial a_1}{\partial x} \right)a_2
+\left(\frac{\partial a_2}{\partial x}- \frac{\partial a_0}{\partial y^\prime} \right)a_1+\left(\frac{\partial a_1}{\partial y^\prime}- \frac{\partial a_2}{\partial y} \right)a_0\neq 0.
\end{equation} 
\end{theorem}
\noindent \textit{Proof.} If such an integrating factor exists, then the conditions  
in Eq. \eqref{muexactconditions} should be hold. A simple calculations shows that the following equations:
\[
a_2\frac{\partial \mu}{\partial y}+\mu\frac{\partial a_2}{\partial y}=a_1\frac{\partial \mu}{\partial y^\prime}+\mu\frac{\partial a_1}{\partial y^\prime}, 
\]

\[
a_2\frac{\partial \mu}{\partial x}+\mu\frac{\partial a_2}{\partial x}=a_0\frac{\partial \mu}{\partial y^\prime}+\mu\frac{\partial a_0}{\partial y^\prime},
\]
and
\[
a_1\frac{\partial \mu}{\partial x}+\mu\frac{\partial a_1}{\partial x}=a_0\frac{\partial \mu}{\partial y}+\mu\frac{\partial a_0}{\partial y},
\]
must be hold. By solving the above three algebraic equations, simultaneously, we get

\[
\left[\left(\frac{\partial a_0}{\partial y}- \frac{\partial a_1}{\partial x} \right)a_2
+\left(\frac{\partial a_2}{\partial x}- \frac{\partial a_0}{\partial y^\prime} \right)a_1+\left(\frac{\partial a_1}{\partial y^\prime}- \frac{\partial a_2}{\partial y} \right)a_0\right]\mu(x,y,z)= 0.
\]
Clearly, if
\[
\left[\left(\frac{\partial a_0}{\partial y}- \frac{\partial a_1}{\partial x} \right)a_2
+\left(\frac{\partial a_2}{\partial x}- \frac{\partial a_0}{\partial y^\prime} \right)a_1+\left(\frac{\partial a_1}{\partial y^\prime}- \frac{\partial a_2}{\partial y} \right)a_0\right]\neq 0,
\]
then
\[
\mu(x,y,y^\prime)=0.
\]
Similarly, for Eq. \eqref{nonlinearsecondorder}, we can show that there is no  integrating factor of one of the forms 
$\mu(x,y)$, $\mu(x,y^\prime)$, or $\mu(y,y^\prime)$ if \eqref{condformu} holds. $\square$

\begin{example}
Consider the second order nonlinear equation
\begin{equation}\label{exmpleq}
xy(2x+y)y^{\prime\prime}+(x^2+xy)y^\prime+(3xy+y^2)=0.
\end{equation}
 Theorem \ref{thm31} shows that the above equation has an integrating factor. In fact, the integrating factor is given by  $\mu(x,y)=\frac{1}{xy(2x+y)}$. This integrating factor transforms Eq. \eqref{exmpleq}  into an exact equation, which can be reduced into a first order differential equation. In fact, it is reduced into the following equation:
\[
\frac{dy}{dx}+\ln\left(xy\sqrt{y+2x}\right)=c.
\]
\end{example}
The following result gives  necessary conditions for the integrating factor to be a function of $x$ only.
\begin{remark}
Through out this paper, we use  the notation $\displaystyle\partial_\eta f:=\frac{\partial f}{\partial \eta}.$
\end{remark}
\begin{lemma} 
 Assume that Eq. \eqref{nonlinearsecondorder} is not an exact equation. Then, it 
has an integrating factor  
\[
\mu(x)=\exp\left\{\int^x\frac{\partial_ya_0-\partial_xa_1}{a_1}dx\right\}=\exp\left\{\int^x\frac{\partial_{y^\prime}a_0-\partial_xa_2}{a_2}dx\right\}
\]
if and only if 
\[
\;\frac{\partial_ya_0-\partial_xa_1}{a_1}\; \text{and}\;\; \;\frac{\partial_{y^\prime}a_0-\partial_xa_2}{a_2}\;
\text{depend only on}\; x,
\]
\[\frac{\partial_ya_0-\partial_xa_1}{a_1}=\frac{\partial_{y^\prime}a_0-\partial_xa_2}{a_2},\]
 and
\[
\;\partial_ya_2=\partial_{y^\prime}a_1. 
\]

\end{lemma} 
 
\noindent\textit{Proof.} Assume that Eq. \eqref{nonlinearsecondorder} has an integrating factor $\mu(x)$. Therefore, conditions  
\eqref{muexactconditions} hold. Hence, we get the following algebraic equations:
\[
\mu\frac{\partial a_2}{\partial y}=\mu\frac{\partial a_1}{\partial y^\prime}, 
\]

\[
a_2\mu^\prime+\mu\frac{\partial a_2}{\partial x}=\mu\frac{\partial a_0}{\partial y^\prime},
\]
and
\[
a_1\mu^\prime+\mu\frac{\partial a_1}{\partial x}=\mu\frac{\partial a_0}{\partial y}.
\] 
Using the first equation, we have a non zero integrating factor, if 
$\frac{\partial a_2}{\partial y}=\frac{\partial a_1}{\partial y^\prime}$. The last two equations implies that
\[
\frac{\mu^\prime}{\mu}=\frac{\frac{\partial a_0}{\partial y^\prime}-\frac{\partial a_2}{\partial x}}{a_2}=\frac{\frac{\partial a_0}{\partial y}-\frac{\partial a_1}{\partial x}}{a_1}.
\]
By integrating the above equation with respect to $x$, we get
\[
\mu(x)=\exp\left\{\int^x\frac{\partial_ya_0-\partial_xa_1}{a_1}dx\right\}=\exp\left\{\int^x\frac{\partial_{y^\prime}a_0-\partial_xa_2}{a_2}dx\right\}. \square
\]  
 
\begin{lemma}
The integrating factor of Eq. \eqref{nonlinearsecondorder} in terms of $y$ is given by  
\[
\mu(y)=\exp\left\{\int^y\frac{\partial_{y^\prime} a_1-\partial_y a_2}{a_2}dy\right\}=\exp\left\{\int^y\frac{\partial_xa_1-\partial_ya_0}{a_0}dy\right\},
\]
provided that 
\[
\;\frac{\partial_{y^\prime} a_1-\partial_y a_2}{a_2}\; \text{and}\;\;\frac{\partial_xa_1-\partial_ya_0}{a_0} \;
\text{depend only on} \;y,
\]
\[
\frac{\partial_{y^\prime} a_1-\partial_y a_2}{a_2}=\frac{\partial_xa_1-\partial_ya_0}{a_0},
\]
and
\[
\;\partial_xa_2=\partial_{y^\prime}a_0. 
\]

\end{lemma} 
 
\begin{lemma}
The integrating factor of Eq. \eqref{nonlinearsecondorder} in terms of $y^\prime$ is given by  
\[
\mu(y^\prime)=\exp\left\{\int^{y\prime}\frac{\partial_{y} a_2-\partial_{y^\prime} a_1}{a_1}dy^\prime\right\}=\exp\left\{\int^{y^\prime}\frac{\partial_xa_2-\partial_{y^\prime}a_0}{a_0}dy^\prime\right\},
\]
provided that 
\[
\;\frac{\partial_{y} a_2-\partial_{y^\prime} a_1}{a_1}\; \text{and}\;\;\frac{\partial_xa_2-\partial_{y^\prime} a_0}{a_0} \;
\text{depend only on}\; y^\prime,
\]
\[
\frac{\partial_{y} a_2-\partial_{y^\prime} a_1}{a_1}=\frac{\partial_xa_2-\partial_{y^\prime}a_0}{a_0},
\]
and
\[
\;\partial_xa_1=\partial_{y}a_0. 
\]
\end{lemma}

\begin{example}
Consider the nonlinear second order differential equation
\[
(1+y^2)yy^{\prime\prime}+g(y)y^\prime+(1+y^2)y=0,
\]
where $g(y)$ is an arbitrary function in $y$. This  equation is not exact. In fact, it has an integrating factor $\mu(y)=\frac{1}{y(1+y^2)}$ which transforms this equation into the exact second order differential equation
\[
y^{\prime\prime}+\frac{g(y)}{y(1+y^2)}y^\prime+1=0.
\]
\end{example}

Since the condition  \eqref{condformu} can not be held easily. i.e., to have an integrating factor of the form $\mu(x,y,y^\prime),$ we are looking for an integrating factor of the form $\mu(\alpha(x)\beta(y)\gamma(y^\prime))$, where $\alpha(x), \beta(y)$ and $\gamma(y^\prime)$ are arbitrary functions in $x,y,$ and $y^\prime$, respectively. For  such an integrating factor to exist, we have the following theorem:
\begin{theorem} 
 Assume that Eq. \eqref{nonlinearsecondorder} is not an exact equation. Then,  an integrating factor  $\mu(\alpha(x)\beta(y)\gamma(y^\prime))$ of Eq. \eqref{nonlinearsecondorder} exists and is given by
\begin{eqnarray*}
\mu(\xi)=\mu(\alpha(x)\beta(y)\gamma(y^\prime))&=&\exp\left\{\int^\xi\frac{\partial_{y^\prime}a_1-\partial_ya_2}{\alpha(x)\left[\beta^\prime(y)\gamma(y^\prime)a_2-\beta(y)\gamma^\prime(y^\prime)a_1\right]}d\xi\right\}\\
&=&\exp\left\{\int^\xi\frac{\partial_ya_0-\partial_xa_1}{\gamma(y^\prime)\left[\alpha(x)\beta^\prime(y)a_1-\alpha^\prime(x)\beta(y)a_0\right]}d\xi\right\}\\
&=&\exp\left\{\int^\xi\frac{\partial_xa_2-\partial_{y^\prime}a_0}{\beta(y)\left[\alpha(x)\gamma^\prime(y^\prime)a_0-\alpha^\prime(x)\gamma(y^\prime)a_2\right]}d\xi\right\},
\end{eqnarray*}
if and only if 
\begin{eqnarray*}
\frac{\partial_{y^\prime}a_1-\partial_ya_2}{\alpha(x)\left[\beta^\prime(y)\gamma(y^\prime)a_2-\beta(y)\gamma^\prime(y^\prime)a_1\right]}
&=&\frac{\partial_ya_0-\partial_xa_1}{\gamma(y^\prime)\left[\alpha(x)\beta^\prime(y)a_1-\alpha^\prime(x)\beta(y)a_0\right]}\\
&=&\frac{\partial_xa_2-\partial_{y^\prime}a_0}{\beta(y)\left[\alpha(x)\gamma^\prime(y^\prime)a_0-\alpha^\prime(x)\gamma(y^\prime)a_2\right]},
\end{eqnarray*}
and they  depend on $\xi(x,y,y^\prime):=\alpha(x)\beta(y)\gamma(y^\prime).$
\end{theorem}
\noindent\textit{Proof.} The proof is a direct consequence of conditions
\eqref{muexactconditions}. 

Using the above theorem, and by either assuming $\gamma(y^\prime)=1$,  $\beta(y)=1$, or $\alpha(x)=1$, we can deduce that the integrating factors are $\mu(\alpha(x)\beta(y))$, $\mu(\alpha(x)\gamma(y^\prime))$ and $\mu(\beta(y)\gamma(y^\prime))$, respectively. The results are listed in the following corollaries:

\begin{corollary} 
 An integrating factor, $\mu(\alpha(x)\beta(y))$, of Eq. \eqref{nonlinearsecondorder} exists and  is given by
\begin{eqnarray*}
\mu(\alpha(x)\beta(y))&=&\exp\left\{\int^\xi\frac{\partial_{y^\prime}a_1-\partial_ya_2}{\alpha(x)\beta^\prime(y)a_2}d\xi\right\}\\
&=&\exp\left\{\int^\xi\frac{\partial_{y^\prime}a_0-\partial_xa_2}{\alpha^\prime(x)
\beta(y)a_2}d\xi\right\}\\&=&\exp\left\{\int^\xi\frac{\partial_ya_0-\partial_xa_1}{\alpha(x)\beta^\prime(y)a_1-\alpha^\prime(x)\beta(y)a_0}d\xi\right\},
\end{eqnarray*}
if and only if 
\[
\;\frac{\partial_{y^\prime}a_1-\partial_ya_2}{\alpha(x)\beta^\prime(y)a_2}=\;
\;\frac{\partial_{y^\prime}a_0-\partial_xa_2}{\alpha^\prime(x)\beta(y)a_2}=\;
\frac{\partial_ya_0-\partial_xa_1}{\alpha(x)\beta^\prime(y)a_1-\alpha^\prime(x)\beta(y)a_0},\;
\]
and they depend on $\xi(x,y):=\alpha(x)\beta(y)$.
\end{corollary}
\begin{corollary}
An integrating factor, $\mu(\alpha(x)\gamma(y^\prime))$, of Eq. \eqref{nonlinearsecondorder} exists and is given by  
\begin{eqnarray*}
\mu(\xi)&=&\mu(\alpha(x)\gamma(y^\prime))\\
&=&\exp\left\{\int^\xi\frac{\partial_y a_2-\partial_{y^\prime}a_1}{\alpha(x)\gamma^\prime(y^\prime)a_0}d\xi\right\}\\
&=&\exp\left\{\int^\xi\frac{\partial_xa_1-\partial_{y}a_0}{\alpha^\prime(x)\gamma(y^\prime)a_1}d\xi\right\}\\
&=&\exp\left\{\int^\xi\frac{\partial_{y^\prime}a_0-\partial_xa_2}{\alpha^\prime(x)\gamma(y^\prime)a_2-\alpha(x)\gamma^\prime(y^\prime)a_0}d\xi\right\},
\end{eqnarray*}
provided that
\[
\frac{\partial_y a_2-\partial_{y^\prime}a_1}{\alpha(x)\gamma^\prime(y^\prime)a_0}=\;
\;\frac{\partial_xa_1-\partial_{y}a_0}{\alpha^\prime(x)\gamma(y^\prime)a_1}=\;
\frac{\partial_{y^\prime}a_0-\partial_xa_2}{\alpha^\prime(x)\gamma(y^\prime)a_2-\alpha(x)\gamma^\prime(y^\prime)a_0},\;
\]
and they  depend on $\xi(x,y^\prime):=\alpha(x)\gamma(y^\prime)$.
\end{corollary} 

\begin{corollary}
An integrating factor, $\mu(\beta(y)\gamma(y^\prime))$, of Eq. \eqref{nonlinearsecondorder} exists and is given by  
\begin{eqnarray*}
\mu(\xi)&=&\mu(\beta(y)\gamma(y^\prime))\\
&=&\exp\left\{\int^\xi\frac{\partial_{y}a_0-\partial_x a_1}{\beta^\prime(y)\gamma(y^\prime)a_1}d\xi\right\}\\
&=&\exp\left\{\int^\xi\frac{\partial_{x}a_2-\partial_{y^\prime}a_0}{\beta(y)\gamma^\prime(y^\prime)a_0}d\xi\right\}\\
&=&\exp\left\{\int^\xi\frac{\partial_{y^\prime}a_1-\partial_ya_2}{\beta^\prime(y)\gamma(y^\prime)a_2-\beta(y)\gamma^\prime(y^\prime)a_1}d\xi\right\},
\end{eqnarray*}
provided that
\[
\frac{\partial_{y}a_0-\partial_x a_1}{\alpha^\prime(y)\beta(y^\prime)a_1}=\frac{\partial_{x}a_2-\partial_{y^\prime}a_0}{\alpha(y)\beta^\prime(y^\prime)a_0}=\frac{\partial_{y^\prime}a_1-\partial_ya_2}{\alpha^\prime(y)\beta(y^\prime)a_2-\alpha(y)\beta^\prime(y^\prime)a_1},
\]
 and they  depend on $\xi(y,y^\prime):=\beta(y)\gamma(y^\prime)$.
\end{corollary} 

\section{Conclusions and Remarks}
In this paper, we imposed conditions on the equation
\[a_2(x,y,y^\prime)y^{\prime\prime}+a_1(x,y,y^\prime)y^\prime+a_0(x,y,y^\prime)=0,\] 
so that it is exact.
  In addition, we introduced an integrating factor in case where the equation is not an exact differential equation. Moreover, we  presented some examples  showing that this method is powerful in solving a class of  second order nonlinear differential equations. For further studies, it is reasonable to improve this definition and this technique to a more complicated class of differential equations. For example, if we consider the general form of the second order nonlinear differential equation $f(x,y,y^\prime,y^{\prime \prime})=0$. Also, it is reasonable to improve this method to work for higher order  nonlinear differential equations.

\end{document}